\documentclass[12pt,a4paper]{article}
\usepackage[T1]{fontenc}
\usepackage{amsmath}
\usepackage{amssymb}
\usepackage{graphicx}
\usepackage{xy}
\title{Free Jordan algebras}
\author{shangshikui}
\date{}
\newcommand{\mZ}{{\mathbb Z}}
\newcommand{\bz}{{\bar{0}}}
\newcommand{\bo}{{\bar{1}}}
\newcommand{\TAG}{{\text{TAG}}}
\newcommand{\TKK}{{\text{TKK}}}
\newcommand{\Lst}{{\text{\bf Liesup}_{\text{T}}}}
\newcommand{\Jsa}{{\text{\bf Jorsup}}}
\newcommand{\slJ}{{sl_2J(D_1|D_2)}}
\newcommand{\ane}{{a_n^{\text{even}}}}
\newcommand{\ano}{{a_n^{\text{odd}}}}

\newcommand{\Hom}{{\text{Hom}}}

\newtheorem{rem}{\noindent\rm\bf Remark}[section]
\newtheorem{cor}{\noindent\rm\bf Corollary}[section]
\newtheorem{thm}{\noindent\rm\bf Theorem}[section]

\newtheorem{pro}{\noindent\rm\bf Proposition}[section]

\newtheorem{lem}{\noindent\rm\bf Lemma}[section]
\newtheorem{conj}{\noindent\rm\bf Conjecture}
\newtheorem{conje}{\noindent\rm\bf Conjecture}

\begin{document}
	
\title{The $\mZ_2$-graded dimensions of the free Jordan superalgebra $J(D_1|D_2)$}

\author{}
\author{Shikui Shang}
%\date{}
	
\maketitle
\vspace{2mm}

\abstract{Let $k$ be a field of characteristic $0$. For a superspace $V=V_\bz\oplus V_\bo$ over $k$, we call the vector $(\dim_k V_\bz ,\dim_k V_\bo )$ the ($\mZ_2$-)graded dimension of $V$. Let $J(D_1|D_2)$ be the free Jordan superalgebra generated by $D_1$ even generators and $D_2$ odd generators. In this paper, we study the graded dimensions of the  $n$-components of $J(D_1|D_2)$ and find the connection between them and the homology of Tits-Allison-Gao Lie superalgebra of $J(D_1|D_2)$ following the method given by I.Kashuba and O.Mathieu in \cite{KM}, where they deal with the free Jordan algebra. And, four interesting conjectures of above contents are proposed in our paper.}

\section{Introduction}

Superspaces and superalgebras appear naturally in theoretical physics and mathematics(\cite{K1} and \cite{Sc}).  Let $k$ be a field of characteristic $0$ and $\mZ_2=\{\bz,\bo\}$. A $\mZ_2$-graded vector 
space $V=V_\bz\oplus V_\bo$ over $k$ is called a superspace. We define its $\mZ_2$-graded dimension by $\overline{\dim }V=(\dim V_\bz ,\dim V_\bo)$ if both $V_\bz$ and $V_\bo$ are finite-dimensional and also call it the graded dimension for simplicity.  Denote the degree or parity of a homogeneous element $v\in V_\alpha$ by $|v|$ for $\alpha\in\{\bz,\bo\}$.

A superalgebra $A$ is a $\mZ_2$-graded algebra $A=A_\bz\oplus A_\bo$. The elements in $A_\bz$ are called even and the elements in $A_\bo$ are called odd. It is well-known that if $V$ is a superspace, then $gl(V)$ becomes a Lie superalgebra.

In the mid-1970s, V.Kac in \cite{K2} introduced the concept of Jordan superalgebras as the generalization of Jordan algebra for super case. Let $J=J_\bz\oplus J_\bo$ be a superalgebra over $k$ with a product $\cdot$. $J$ is call a Jordan superalgebra if it satisfies the following two identities for any homogeneous $x,y,z\in J$
\begin{align*}
(1)&(\text{super commutativity}) x\cdot y-(-1)^{|x||y|}y\cdot x=0,\\
(2)&(\text{super Jordan identity)}\\
&(-1)^{|x||z|}[L_{x\cdot y},L_z]+(-1)^{|y||x|}[L_{y\cdot z},L_x]+(-1)^{|z||y|}[L_{z\cdot x},L_y]=0
\end{align*}
where $L_x$ is an endomorphism of $J$ with $L_x(y)=x\cdot y$ and the bracket is taken in the Lie superalgebra $gl(J)$. Some theories of simple and semisimple finite dimensional Jordan superalgebras are studied in \cite{CK}, \cite{Ka2}, \cite{MZ},  and \cite{Z2}. 

Let $J(D_1|D_2)$ be the free Jordan superalgebra (without unit) over $k$ generated by $D_1$ even generators $x_1,\cdots,x_{D_1}$ and $D_2$ odd generators $y_1,\cdots,y_{D_2}$. Then, $J(D_1|D_2)$ has a natural $\mZ_{>0}$-grading 
$$J(D_1|D_2)=\oplus_{n\geq 1}J_n(D_1|D_2).$$
For any $n\geq 1$, $J_n(D_1|D_2)$ is a superspace with a finite graded dimension. When $D_2=0$, the free Jordan algebra $J(D_1)$ has been investigated in many papers by the Novosibirsk school of algebra, e.g. \cite{Me}, \cite{Z}, \cite{ZSSS}. 

Recently, the dimension of the $n$-component $J_n(D_1)$ of $J(D_1)$ is studied using Lie algebra homology theory by I. Kashuba and O. Mathieu in paper \cite{KM}. Their work inspires us to consider the free Jordan superalgebra case. For $n=1$, the graded dimension of $J_1(D_1|D_2)$ is $(D_1,D_2)$. But for $n\geq 2$, the vector $\overline{\dim}J_n(D_1|D_2)$ is difficult to deal. We will give a conjecture(Conjecture 4 in Section 4 or 6) on $J(D_1|D_2)$ which provides a recurrence relation $(\mathcal E)$ to uniquely determine $\overline{\dim}J_n(D_1|D_2)$. 

The theory of Jordan algebras is very closed associative with the theory of $sl_2$-graded($3$-graded) Lie algebras(see \cite{J},\cite{Ka},\cite{Ko} and \cite{T}). There are also similar relations between Jordan superalgebras and $sl_2$-graded Lie superalgebras(see \cite{BC},\cite{KiM} and \cite{KMZ}). The conjecture on the graded dimension of $J_n(D_1|D_2)$ can be deduced from a conjecture(Conjecture 1 in Section 3 or 6) about Lie superalgebra homology.

In detail, for $3$-dimensional Lie algebra
$sl_2(k)$ and $m\geq 0$, let $L(m)$ be the simple module of $sl_2(k)$ with dimension $m+1$. For an $sl_2$-module 
$M$, denote $M^{sl_2}$ the sum of all submodules which are isomorphic to $L(0)$ and   $M^{\text{ad}}$ the sum of all submodules which are isomorphic to $L(2)$. We consider the category $\Lst$ of Lie superalgebras ${\mathfrak g}$ over $k$ which are also $sl_2(k)$-modules, on which the $sl_2(k)$ acts as (even) derivations such that ${\mathfrak g}$ has a decomposition ${\mathfrak g}={\mathfrak g}^{sl_2}\oplus {\mathfrak g}^{\text{ad}}$. Let $\Jsa$ be the category of Jordan superalgebras over $k$. There are two functors, Tits-Allison-Gao functor $\TAG$ and Tits functor $T$, constituting an adjoint pair of functors between categories $\Jsa$ and $\Lst$.

We denote $\TAG(J(D_1|D_2))$ simply by $\slJ$. Since $J(D_1|D_2)$ is free and the $\TAG$-construction is functorial, it is very natural to expect some homology superspaces of $\TAG(J(D_1|D_2))$ are vanishing. This is Conjecture 1 in Section 3. We will see that the Conjecture 1 is stronger which can implies Conjecture 4.

The paper is organized as follows. In section 2, we give the details of the adjoint pair of Tits-Allison-Gao functor $\TAG$ and Tits functor $T$ for the super case. The homology theory of Lie superalgebra is introduced in Section 3. Moreover, $H_r(\slJ)$ for $r=0,1,2$ are computed and a conjecture(Conjecture 1) on the higher homology superspaces is stated. In section 4, more structures on $J(D_1|D_2)$ are studied and  three conjectures(Conjecture 2-4) are proposed. We show the relations between the homology of $\slJ$ and them. In section 5, some numerical evidences of above conjectures are given. Finally, we summarize the main results in the last section of this paper.

\section{The Tits-Allison-Gao functor and Tits functor for super case}
\subsection{The TKK-construction}
Let $J=J_\bz\oplus J_\bo$ be a Jordan superalgebra over $k$ which is not necessarily unital. 
$gl(J)$ is the Lie superalgebra of endomorphisms on $J$. For $\alpha\in\{\bz,\bo\}$, $D\in gl(J)_\alpha$ is called a (super-)derivation of $J$ of degree $|D|=\alpha$, if 
$$D(xy)=D(x)\cdot y+(-1)^{|D||x|}x\cdot D(y)$$ 
for any homogeneous $x,y\in J$ and denote $D\in\text{Der}(J)_\alpha$. Set the Lie superalgebra
$$\text{Der}(J)=\text{Der}(J)_\bz\oplus \text{Der}(J)_\bo.$$

By the definition of Jordan superalgebras, $$\partial_{x,y}:=[L_x,L_y]=L_xL_y-(-1)^{|x||y|}L_yL_x\in\text{Der}(J)_{|x|+|y|}.$$ The linear combinations of $\partial_{x,y}$ are called the inner derivations of $J$, and denote the subspace $\text{Inn}(J)$ of inner derivations. Since for any homogeneous $D\in\text{Der}(J)$ and  $x,y\in J$, 
$$[D,\partial_{x,y}]=\partial_{D(x),y}+(-1)^{|D||x|}\partial_{x,D(y)},$$ $\text{Inn}(J)$ is a graded ideal of $\text{Der}(J)$.

\ 

{\bf Tits-Kantor-Koecher Construction:}
The Tits-Kantor-Koecher algebra $\TKK(J)$ of a Jordan superalgebra $J$ is defined as follows. As a vector space
$$\TKK(J)=\text{Inn}(J)\oplus(sl_2(k)\otimes J)$$
with $\mZ_2$-grading
$$\TKK(J)_\alpha=\left\{\sum_i\partial_{x_i,y_i} \ |\ |x_i|+|y_i|=\alpha\right\}\oplus(sl_2(k)\otimes J_\alpha)$$
for $\alpha\in\mZ_2$. Let $\kappa:sl_2(k)\times sl_2(k)\rightarrow k$ be the Killing form on $sl_2(k)$. The bracket of $\TKK(J)$ is given such that $\text{Inn}(J)$ becomes a subalgebra of $\TKK(J)$ and for $a,b\in sl_2(k), x,y,z\in J, D\in\text{Inn}(J)$ 
\begin{align*}
&1.\ \ \ \ [a\otimes x,b\otimes y]=\frac{1}{2}\kappa(a,b)\partial_{x,y}+[a,b]\otimes(x\cdot y),\\	
&2.\ \ \ \ [D, a\otimes z]=a\otimes D(z).
\end{align*}

It is proved in \S4.4 of \cite{BC} that $\TKK(J)$ is a Lie superalgebra over $k$. However, $\TKK$ is not functorial since the notion of inner
derivations is not functorial.

\subsection{The Lie superalgebra $\TAG(J)$}

Following \cite{KM}, we introduce the Tits-Allison-Gao construction of a Jordan superalgebra $J$, which provides a refinement of the $\TKK$-construction and was first studied in the paper \cite{AG}.

Let $J=J_\bz\oplus J_\bo$ be a Jordan superalgebra. Then, $J\otimes_kJ$ is a superspace with the canonical $\mZ_2$-grading. Set $R^s(J)$ is the subspace of $J\otimes_kJ$ spanned by all homogeneous elements of forms
\begin{align*}&I_{x,y}=x\otimes y+(-1)^{|x||y|}y\otimes x,\\
&I_{x,y,z}=(-1)^{|x||z|}(x\cdot y)\otimes z+(-1)^{|y||x|}(y\cdot z)\otimes x+(-1)^{|z||y|}(z\cdot x)\otimes y,
\end{align*}
where $x,y,z$ run over the homogeneous elements in $J$. And, set the superspace $${\mathcal B}^s(J)=(J\otimes_kJ)/R^s(J),$$
and denote the canonical image of $x\otimes y$ in ${\mathcal B}^s(J)$ by $\{x\otimes y\}$.

\begin {lem}\label{lem:2.1}The $k$-linear map
$$\varphi:{\mathcal B}^s(J)\rightarrow\text{Inn}(J),\{x\otimes y\}\mapsto\partial_{x,y}$$
are surjective.
\end {lem}
{\bf Proof. }Since $\partial_{x,y}$ are bilinear in $x,y$ and $\text{Inn}(J)$ is spanned by $\partial_{x,y}$,
$J\otimes J\rightarrow\text{Inn}(J), x\otimes y\mapsto\partial_{x,y}$ are well-defined and surjective.

By the definition of Jordan superalgebras, $R^s(J)$ is contained in the kernel of the above map. Hence, we have the surjective map
$\varphi:{\mathcal B}^s(J)\rightarrow\text{Inn}(J),\{x\otimes y\}\mapsto\partial_{x,y}$. $\Box$  

The superspace $J\otimes_kJ$ is a graded $\text{Inn}(J)$-module induced by the natural action of $\text{Inn}(J)$ on $J$, i.e.,
$$\partial_{x,y}(z\otimes w)=\partial_{x,y}(z)\otimes w+(-1)^{(|x|+|y|)|z|}z\otimes\partial_{x,y}(w).$$

\begin {lem}\label{lem:2.2}$R^s(J)$ is a graded submodule of $J\otimes_kJ$ and the quotient ${\mathcal B}^s(J)$ is also an $\text{Inn}(J)$-module.
\end {lem}
{\bf Proof. }Using direct computation, we have 
\begin{align*}
&\partial_{x,y}(I_{z,w})=I_{\partial_{x,y}(z),w}+(-1)^{(|x|+|y|)|z|}I_{z,\partial_{x,y}(w)},\\
&\partial_{x,y}(I_{u,v,w})=\\
&(-1)^{(|x|+|y|)|w|}I_{\partial_{x,y}(u),v,w}+(-1)^{(|x|+|y|)|u|}I_{u,\partial_{x,y}(v),w}+(-1)^{(|x|+|y|)|v|}I_{u,v,\partial_{x,y}(w)},
\end{align*} 
for homogeneous $x,y,z,u,v,w\in J$.

Hence, $R^s(J)$ is invariant under the action of $\text{Inn}(J)$, and ${\mathcal B}^s(J)=(J\otimes_k J)/R^s(J)$ is the quotient $\text{Inn}(J)$-module. $\Box$

We define a bracket product on ${\mathcal B}^s(J)$ by
$$[\{x\otimes y\},\{z\otimes w\}]=\partial_{x,y}(\{z\otimes w\}).$$

We claim that ${\mathcal B}^s(J)$ is a Lie superalgebra over $k$ with such bracket. The following lemma is needed,

\begin {lem}\label{lem:2.3}Let $M$ be a graded module over a Lie superalgebra $\mathfrak g$ and $\lambda:M\rightarrow{\mathfrak g}$ a $\mathfrak g$-module homomorphism. Then,
$$A(M)=\text{span}\{\lambda(m).n+(-1)^{|m||n|}\lambda(n).m\ |\ m,n \text{ homogeneous in }M\}$$
is a submodule of $M$ satisfying $\text{Im}(\lambda).\ker(\lambda)\subseteq A(M)\subseteq\ker(\lambda)$. Moreover, the quotient module $Q=M/A(M)$ forms a Lie superalgebra relative to the product $[p,q]=\lambda(p).q$ for $p,q\in Q$ and the map $\mu:M/A(M)\rightarrow\mathfrak g$ induced by $\lambda$ is a central extension of $\mathfrak g$.
\end {lem}
{\bf Proof. }The proof is parallel to the Lie algebra version of Lemma 3.1 in \cite{S}. $\Box$

\begin {cor}\label{cor:2.1} ${\mathcal B}^s(J)$ is a Lie superalgebra over $k$ and $\varphi:{\mathcal B}^s(J)\rightarrow\text{Inn}(J)$ is a central extension of Lie superalgebras.
\end {cor}
{\bf Proof. }Applied Lemma \ref{lem:2.2} and Lemma \ref{lem:2.3}, we only need to show that $A({\mathcal B}^s(J))=0$.
Its verification is straightforward as follows. For homogeneous $x,y,z,w\in J$, in $J\otimes_kJ$,
\begin{align*}
&\partial_{x,y}(z\otimes w)+(-1)^{(|x|+|y|)(|z|+|w|)}\partial_{z,w}(x\otimes y)\\
=&x\cdot(y\cdot z)\otimes w-(-1)^{|x||y|}y\cdot(x\cdot z)\otimes w+\\
&(-1)^{(|x|+|y|)|z|}z\otimes x\cdot(y\cdot w)-(-1)^{(|x|+|y|)|z|+|x||y|}z\otimes y\cdot(x\cdot w)+\\
&(-1)^{(|x|+|y|)(|z|+|w|)}z\cdot(w\cdot x)\otimes y-(-1)^{(|x|+|y|)(|z|+|w|)+|z||w|}w\cdot(z\cdot x)\otimes y\\
+&(-1)^{|y|(|z|+|w|)}x\otimes z\cdot(w\cdot y)-(-1)^{|y||z|+|y||w|+|z||w|}x\otimes w\cdot(z\cdot y)\\
\equiv
%&(-1)^{|x||w|}I_{x,y\cdot z,w}-(-1)^{(|x|+|w|)|y|}I_{y,x\cdot z,w}+(-1)^{|x|(|z|+|w|)+|y||w|}I_{z,w\cdot x,y}\\
&-(-1)^{(|x|+|y|+|z|)|w|}(w\cdot x)\otimes(y\cdot z)+(-1)^{(|x|+|y|+|z|)|w|+|x||y|}(w\cdot y)\otimes(x\cdot z)\\
&-(-1)^{|x|(|y|+|z|+|w|)+|z||w|}(y\cdot w)\otimes(z\cdot x)+(-1)^{(|y|+|z|)|w|+|y||z|}(x\cdot w)\otimes(z\cdot y)\\
\equiv&0(\mod R^s(J)). \ \ \ \ \Box
\end{align*}

\

{\bf Tits-Allison-Gao Construction:} For a Jordan superalgebra $J$,
the Tits-Allison-Gao algebra $\TAG(J)$ is defined as the a vector space
$$\TAG(J)={\mathcal B}^s(J)\oplus(sl_2(k)\otimes J)$$
with $\mZ_2$-grading
$$\TAG(J)_\alpha=\left\{\sum_i\{x_i\otimes y_i\} \ |\ |x_i|+|y_i|=\alpha\right\}\oplus(sl_2(k)\otimes J_\alpha)$$
for $\alpha\in\mZ_2$ and the bracket product
\begin{align*}
	&1.\ \ \ \ [a\otimes x,b\otimes y]=\frac{1}{2}\kappa(a,b)\{x\otimes y\}+[a,b]\otimes(x\cdot y),\\	
	&2.\ \ \ \ [\{x\otimes y\}, a\otimes z]=a\otimes\partial_{x,y}(z),\\
	&3.\ \ \ \ [\{x\otimes y\},\{z\otimes w\}]=\partial_{x,y}(\{z\otimes w\}),
\end{align*}
 for $a,b\in sl_2(k), x,y,z,w\in J$. 

\begin{thm}\label{thm:2.1} $\TAG(J)$ is a Lie superalgebra over $k$ for Jordan superalgebra $J$.
\end{thm}
{\bf Proof. }In Corollary \ref{cor:2.1}, we have seen that ${\mathcal B}^s(J)$ is a Lie superalgebra, which is a subalgebra of $\TAG(J)$.
It is routine to check the rest cases for elements in $\TAG(J)$ satisfy the super commutativity and super Jordan identity. $\Box$

\subsection{The Category $\Lst$ and Tits functor $T$}

Let $\text{\bf T}$ be the category of $sl_2(k)$-modules $M$ such that $M=M^{sl_2}\oplus M^{\text{ad}}$. In this paper, weights of an $sk_2(k)$-module $M$ always mean the eigenvalues of $h\in sl_2(k)$.

 Let $\Lst$ be the category of Lie superalgebras ${\mathfrak g}={\mathfrak g}_\bz\oplus{\mathfrak g}_\bo$ in $\text{\bf T}$ on which $sl_2(k)$ acts as even derivations. Particularly,  both ${\mathfrak g}_\bz$ and ${\mathfrak g}_\bo$ are $sl_2(k)$-submodules.

Let $J$ be a Jordan super algebra over $k$. The algebra $\TKK(J)$ and $\TAG(J)$ are objects in $\Lst$, where
$$\TKK(J)^{sl_2}=\text{Inn}(J),\ \ \ \TKK(J)^{\text{ad}}=sl_2(k)\otimes J$$
and
$$\TAG(J)^{sl_2}={\mathcal B}^s(J), \ \ \ \TAG(J)^{\text{ad}}=sl_2(k)\otimes J.$$

Let $\{e,h,f\}$ be the usual basis o $sl_2(k)$. For a Lie superalgebra $\mathfrak g$ in $\Lst$, set
$$T({\mathfrak g})=\{x\in{\mathfrak g}\ |\ h.x=2x\}.$$
We define a new product $\cdot$ between two elements $x,y\in T({\mathfrak g})$ by
$$x\cdot y=\frac{1}{2}[x,f.y]\in T({\mathfrak g}).$$
\begin{thm}\label{thm:2.2} $T({\mathfrak g})$ is a Jordan superalgebra over $k$ with above dot product and $T$ is a functor from the category $\Lst$ to the category $\Jsa$ of Jordan superalgebra over $k$.
\end{thm}
{\bf Proof. }We only need to show $T({\mathfrak g})$ is a Jordan superalgebra since it is obvious that $T$ is functorial. It follows the proof of Proposition 1 in \cite{T}  for Tits functor on Lie algebras.

Since the weights of $\mathfrak g$ are not greater than $2$ and $sl_2(k)$ acts as even derivations on $\mathfrak g$, for $x,y\in T({\mathfrak g})$, we have that $[x,y]=0$ and 
$$x\cdot y-(-1)^{|x||y|}y\cdot x=\frac{1}{2}([x,f.y]+[f.x,y])=f.[x,y]=0.$$
Hence, the super commutativity holds.

Next, we check the super-Jordan identity. For $x\in J$, using the notations
$$e(x)=x,\ \ h(x)=-f.e(x),\ \ f(x)=\frac{1}{2}f.h(x),$$
we have that $a.b(x)=[a,b](x)$ and $|a(x)|=|x|$ for any $a,b\in sl_2(k)$ and $x\in J$.

Then, for any $x,y\in J$,
\begin{equation}[h(x),e(y)]=-[f.e(x),e(y)]=[e(x),f.e(y)]=2e(x\cdot y)\tag{2.1}\end{equation} and 
\begin{align*}&[e(x),f(y)]=\frac{1}{2}[e(x), f.h(y)]=\frac{1}{2}\left(f.[e(x),h(y)]-[f.e(x),h(y)]\right)\\
	=&\frac{1}{2}\left(-f.[e(x),f.e(y)]-[f.e(x),h(y)]\right)=\frac{1}{2}\left(-2f.e(x\cdot y)+[h(x),h(y)]\right)\\
	=&h(x\cdot y)+\frac{1}{2}[h(x),h(y)].\tag{2.2}
\end{align*}
Moreover, $e.[e(x),f(y)]=-2e(x\cdot y)+\frac{1}{2}e.[h(x),h(y)]$. 

On the other hand, we see that
$$e.[e(x),f(y)]=[e(x),e.f(y)]=[e(x),h(y)]=-2e(x\cdot y).$$
Comparing with above equation, one has that $e.[h(x),h(y)]=0$ and $[h(J),h(J)]\subseteq{\mathfrak g}^{sl_2}$.
Then,
\begin{align*}&[h(x),f(y)]=-[f.e(x),f(y)]=-f.[e(x),f(y)]\\
	=&-f.(h(x\cdot y)+\frac{1}{2}[h(x),h(y)])=-f.h(x\cdot y)\\
	=&-2f(x\cdot y).\tag{2.3}
\end{align*}

For $x,y,z\in J$, by Equation (2.2) and super-commutativity of $\cdot$,
\begin{align*}
	&[[e(x),f(y)],h(z)]=[e(x),[f(y),h(z)]]+(-1)^{|y||z|}[[e(x),h(z)],f(y)]\\
	=&2([e(x),f(y\cdot z)]-(-1)^{|y||z|}[e(x\cdot z),f(y)])\\
	=&2h(x\cdot(y\cdot z))-(-1)^{|y||z|}(x\cdot z)\cdot y))\\
	&+[h(x),h(y\cdot z)]-(-1)^{|y||z|}[h(x\cdot z),h(y)].
\end{align*}
Hence, by Equation (2.2),
\begin{align*}
	&[h(x\cdot y),h(z)]=\left[[e(x),f(y)]-\frac{1}{2}[h(x),h(y)],h(z)\right]\\
	=&[[e(x),f(y)],h(z)]-\frac{1}{2}[[h(x),h(y)],h(z)]\\
	=&2h(x\cdot(y\cdot z))-(-1)^{|y||z|}(x\cdot z)\cdot y))\\
	&+[h(x),h(y\cdot z)]-(-1)^{|y||z|}[h(x\cdot z),h(y)]-\frac{1}{2}[[h(x),h(y)],h(z)]. \tag{2.4}
\end{align*}

Acting the even permutations in the cyclic group $C_3$ on $\{x,y,z\}$ in Equation (2.4) and using the super commutativity, we obtain that
\begin{align*}
	&(-1)^{|x||z|}[h(x\cdot y),h(z)]+(-1)^{|y||x|}[h(y\cdot z),h(x)]+(-1)^{|z||y|}[h(z\cdot x),h(y)]\\
	=&2h((-1)^{|x||z|}x\cdot(y\cdot z)-(-1)^{(|x|+|y|)|z|}(x\cdot z)\cdot y\\
	&+(-1)^{|y||x|}y\cdot(z\cdot x)-(-1)^{(|y|+|z|)|x|}(y\cdot x)\cdot z\\
	&+(-1)^{|z||y|}z\cdot(x\cdot y)-(-1)^{|(z|+|x|)|y|}(z\cdot y)\cdot x)\\
	&+(-1)^{|x||z|}[h(x),h(y\cdot z)]-(-1)^{(|x|+|y|)|z|}[h(x\cdot z),h(y)]\\
	&+(-1)^{|y||x|}[h(y),h(z\cdot x)]-(-1)^{(|y|+|z|)|x|}[h(y\cdot x),h(z)]\\
	&+(-1)^{|z||y|}[h(z),h(x\cdot y)]-(-1)^{(|z|+|x|)|y|}[h(z\cdot y),h(x)].
\end{align*}
which is equivalent that
$$3((-1)^{|x||z|}[h(x\cdot y),h(z)]+(-1)^{|y||x|}[h(y\cdot z),h(x)]+(-1)^{|z||y|}[h(z\cdot x),h(y)])=0.$$

Since the characteristic of $k$ is zero, $$(-1)^{|x||z|}[h(x\cdot y),h(z)]+(-1)^{|y||x|}[h(y\cdot z),h(x)]+(-1)^{|z||y|}[h(z\cdot x),h(y)]=0.$$

Meanwhile, by Equation (2.1),  
\begin{align*}&[[h(x),h(y)],e(z)]=4e(x\cdot(y\cdot z)-(-1)^{|x||y|}y\cdot(x\cdot z))\\
	=&4e(L_xL_y(z)-(-1)^{|x||y|}L_yL_x(z))=4e([L_x,L_y](z)).\end{align*}

As operators on $T({\mathfrak g})$, we have 
$$(-1)^{|x||z|}[L_{x\cdot y},L_z]+(-1)^{|y||x|}[L_{y\cdot z},L_x]+(-1)^{|z||y|}[L_{z\cdot x},L_y]=0.$$ 
Therefore, the super Jordan identity holds and $T({\mathfrak g})$ is a super Jordan aglebra. $\Box$

\begin {lem}\label{lem:2.4}Let ${\mathfrak g}\in\Lst$. Then, there is a Lie superalgebra morphism
$$\theta_{\mathfrak g}:\TAG(T({\mathfrak g}))\rightarrow{\mathfrak g}$$
which is the identity on $T(\mathfrak g)$.
\end {lem}

{\bf Proof. }
First, we have
$${\mathfrak g}={\mathfrak g}^{sl_2}\oplus sl_2\otimes T({\mathfrak g}).$$
Since $\Hom_{sl_2}(sl_2^{\otimes 2}, k)=k.\kappa$, there is a linear map $\psi:T({\mathfrak g})\otimes T({\mathfrak g})\rightarrow{\mathfrak g}^{sl_2}$ such that 
$$[a(x),b(y)]=[a,b](x\cdot y)+\frac{\kappa(a,b)}{2}\psi(x,y)$$
for $a,b\in sl_2(k)$ and $x,y\in T({\mathfrak g})$.

Since $[h(x),h(y)]+(-1)^{|x||y|}[h(y),h(x)]=0$ and $\kappa(h,h)=8\neq0$,
\begin{equation*}\psi(x,y)+(-1)^{|x||y|}\psi(y,x)=0\tag{2.5}\end{equation*}
for $x,y\in T({\mathfrak g})$.

For $a,b,c\in sl_2(k)$, we have that
$$[[a(x),b(y)],c(z)]=[[a,b],c]((x\cdot y)\cdot z)+\kappa([a,b],c)\psi(x\cdot y,z)+\kappa(a,b)[\psi(x,y),c(z)].$$

Using the super Jacobi identity for triple $(e(x),f(y),h(z))$ in $\mathfrak g$, \begin{align*}&(-1)^{|x||z|}[[e(x),f(y)],h(z)]+(-1)^{|y||x|}[[f(y),h(z)],e(x)]+\\
	&(-1)^{|z||y|}[[h(z),e(x)],f(y)]=0.
\end{align*}
Since $\kappa(h,h)=2\kappa(e,f)=8\neq 0$, the component of ${\mathfrak g}^{sl_2}$ in above equation gives us for $x,y,z\in T({\mathfrak g})$,
\begin{equation*}(-1)^{|x||z|}\psi(x\cdot y,z)+(-1)^{|y||x|}\psi(y\cdot z,x)+(-1)^{|z||y|}\psi(z\cdot x,y)=0.\tag{2.6}\end{equation*}

By (2.5) and (2.6), we have that the map $\psi$ factors through ${\mathcal B}^s(T({\mathfrak g}))$. A linear map $\theta_{\mathfrak g}:\TAG(T({\mathfrak g}))\rightarrow{\mathfrak g}$ is defined by requiring that $\theta_{\mathfrak g}$ is identity on $sl_2\otimes T({\mathfrak g})$ and $\theta_{\mathfrak g}=\psi$ on ${\mathcal B}^s(T({\mathfrak g}))$. It is easy to check that $\theta_{\mathfrak g}$
is a morphism of Lie superalgebra over $k$.   $\Box$ 

It is clear that $\TAG:\Jsa\rightarrow\Lst, J\mapsto\TAG(J)$ is also a functor. More precisely, we have 
\begin{thm}\label{thm:2.3} The functor $\TAG:\Jsa\rightarrow\Lst$ is the left adjoint of the Tits functor $T$, namely,
$$\Hom_\Lst(\TAG(J),{\mathfrak g})=\Hom_\Jsa(J,T({\mathfrak g}))$$
for any $J\in\Jsa$ and ${\mathfrak g}\in\Lst$.
\end{thm}
{\bf Proof. } Since $T(\TAG(J))=J$, there is a natural map by restricting on $J$
$$\mu:\Hom_{\Lst}(\TAG(J), {\mathfrak g})\rightarrow\Hom_{\Jsa}(J, T({\mathfrak g})), \eta\mapsto\eta|_J.$$
Since $\TAG(J)$ is generated by $sl_2\otimes J$, $\mu$ is injective. 

On the other hand, let $\phi:J\rightarrow T(\mathfrak g)$
be a morphism of Jordan superalgebra. By functoriality of the $\TAG$-construction, we get a Lie superalgebra morphism
$$\TAG(\phi):\TAG(J)\rightarrow\TAG(T({\mathfrak g}))$$
and by Lemma \ref{lem:2.4} there is a canonical Lie superalgebra morphism
$$\theta_{\mathfrak g}:\TAG(T({\mathfrak g}))\rightarrow{\mathfrak g}.$$

So, $\theta_{\mathfrak g}\circ\TAG(\phi)$ extends $\phi$ to a morphism of Lie superalgebras. Therefore, $\phi=\mu(\theta_{\mathfrak g}\circ\TAG(\phi))$ and $\mu$ is bijective.
$\Box$

\section{The homology of $\TAG$ algebra $\slJ$}

\subsection{The (co)homology of Lie superalgebras in general}

The (co)homology theory of Lie superalgebras over $k$ is parallel to the associated theory of Lie algebras(See \cite{AR-G} \cite{BB} and \cite{SZ}). 

For a Lie superalgebra $\mathfrak g$, the functor $-_{\mathfrak g}$(resp. $-^{\mathfrak g}$) from the category of $\mathfrak g$-module $M$ to the category of superspaces over $k$ is given by
$$M_{\mathfrak g}=M/{\mathfrak g}M, \ \ (\text{resp. } M^{\mathfrak g}=\{m\in M\ |\ x.m=0, \forall x\in{\mathfrak g}\})$$

The functor $-_{\mathfrak g}$(resp. $-^{\mathfrak g}$) is right(resp. left) exact. For $r\geq 0$, the $r$-th right(resp. left) derived functor $H_r({\mathfrak g},-)$(resp. $H^r({\mathfrak g},-)$) of  
 $-_{\mathfrak g}$(resp. $-^{\mathfrak g}$) is the called the $r$-th homology(resp. cohomology) functor.
 And, $H_r({\mathfrak g},M)$(resp. $H^r({\mathfrak g},M)$) is called the $r$-th homology(resp. cohomology) superspace of $\mathfrak g$ with coefficients in $M$.
 
In particular, we set $H_r({\mathfrak g})=H_r({\mathfrak g},k)$(resp. $H^r({\mathfrak g})=H^r({\mathfrak g},k)$) and call them the $r$-th homology(resp. cohomology) superspaces of $\mathfrak g$ where $k$ is viewed as the trivial $1$-dimensional module of $\mathfrak g$.

By the property of derived functors, 
\begin {lem}\label{lem:3.1} For any Lie superalgebra $\mathfrak g$,
$$H_0({\mathfrak g})=H^0({\mathfrak g})=k. \ \ \ \Box$$
\end {lem}

The Chevalley-Eilenberg complex can be used to compute the (co)homology of Lie superalgebra ${\mathfrak g}={\mathfrak g}_\bz\oplus{\mathfrak g}_\bo$. Let $U(\mathfrak g)$ be the enveloping algebra of $\mathfrak g$. For $r\geq 0$, set superspace
$$V_r({\mathfrak g})=\oplus_{p+q=r}\Lambda^p({\mathfrak g}_\bz)\oplus S^q({\mathfrak g}_\bo).$$
with $V_r({\mathfrak g})_\alpha=\oplus_{\stackrel{p+q=r}{q\equiv\alpha(\mod 2)}}\Lambda^p({\mathfrak g}_\bz)\otimes S^q({\mathfrak g}_\bo)$ for $\alpha\in\mZ_2$ and the differential $d_r:V_r({\mathfrak g})\rightarrow V_{r-1}({\mathfrak g})$ is given by
\begin{align*}
&d_r((x_1\wedge\cdots\wedge x_p)\otimes(y_1\cdots y_q))\\
=&\sum_{1\leq i<j\leq p}(-1)^{i+j+1}([x_i,x_j]\wedge x_1\wedge\cdots\wedge\hat{x}_i\wedge\cdots\wedge\hat{x}_j\wedge\cdots x_p)\otimes(y_1\cdots y_q)\\
&+\sum_{i=1}^p\sum_{j=1}^q(-1)^{i+1}(x_1\wedge\cdots\wedge\hat{x}_i\wedge\cdots\wedge x_{p-1})\otimes([x_i,y_j]y_1\cdots \hat{y}_j\cdots y_q)\\
&+\sum_{1\leq i<j\leq q}([y_i,y_j]\wedge x_1\wedge\cdots\wedge x_p)\otimes(y_1\cdots \hat{y}_j\hat{y}_{j+1}\cdots y_q).
\end{align*}
Then, for a graded ${\mathfrak g}$-module $M$,
\begin{align*}
	&H_r({\mathfrak g},M)=\frac{\ker(V_r({\mathfrak g})\otimes M\rightarrow V_{r-1}({\mathfrak g})\otimes M)}{\text{Im}((V_{r+1}({\mathfrak g})\otimes M\rightarrow V_r({\mathfrak g})\otimes M)},\\
	&H^r({\mathfrak g},M)=\frac{\ker(\Hom_k(V_r({\mathfrak g}),M)\rightarrow\Hom_k(V_{r+1}({\mathfrak g}),M)}{\text{Im}(\Hom_k(V_{r-1}({\mathfrak g}),M)\rightarrow\Hom_k(V_r({\mathfrak g}),M)}.
\end{align*}
where all morphisms are induced by the above differentials $d_*$.

Let ${\mathfrak g}\in\Lst$ and $M$ be a ${\mathfrak g}$-module with compatible $sl_2(k)$-actions. Then,  $H_r({\mathfrak g},M)$ and $H^r({\mathfrak g},M)$
are also $sl_2(k)$-modules. Note that even for ${\mathfrak g}\in\Lst$, the weights of (co)homology superspaces may be greater than $2$.

For Lie superalgebra $\TAG(J)$, one has 
\begin {lem}\label{lem:3.2} For any Jordan superalgebra $J$,
$$H_1({\TAG(J)})=sl_2(k)\otimes(J/J^2).$$
Particularly, $H_1({\TAG(J)})=H_1({\TAG(J)})^{\text{ad}}$.
\end {lem}
{\bf Proof.} Using the Chevalley-Eilenberg complex, we have $H_1({\mathfrak g})={\mathfrak g}/[{\mathfrak g},{\mathfrak g}]$. Taking ${\mathfrak g}=\TAG(J)$, 
it is easy to compute
$$[\TAG(J),\TAG(J)]={\mathcal B}^s(J)\oplus(sl_2(k)\otimes J^2).$$
Hence, $H_1({\TAG(J)})=sl_2(k)\otimes(J/J^2)$.    $\Box$

\subsection{The homology of $\slJ$}

Next, we consider the free Jordan superalgebra $J=J(D_1|D_2)$ and denote $\slJ=\TAG(J(D_1|D_2))$ for simplicity. Endowing every generator degree $1$,
$J(D_1|D_2)$ has a natural $\mZ_+$-grading, namely,
$$J(D_1|D_2)=\oplus_{n\geq 1}J_n(D_1|D_2).$$

To compute the $H_2(\slJ)$, we need to consider the central extension of $\slJ$. The following lemma is the super version of  Lemma 13 in
\cite{KM}.
\begin {lem}\label{lem:3.3} Let $M$ be a finite-dimensional $\slJ$-module. Assume that the $sl_2(k)$-actions on $M$ are compatible with the $\slJ$-actions, i.e., $a.(\rho.m)=(a.\rho).m+\rho.(a.m)$ for all $a\in sl_2(k),\rho\in\slJ, m\in M$. Then,
$$H_2(\slJ,M)^{sl_2}=0.$$\\
\end {lem}
{\bf Proof.} By duality, we only need to show $H^2(\slJ,M^*)^{sl_2}=0$.

A super $2$-cocycle $\omega:\slJ\times\slJ\rightarrow M^*$ is called $sl_2(k)$-invariant if 
$$a.\omega(\rho,\sigma)=\omega(a.\rho,\sigma)+(\rho,a.\sigma)$$
for $a\in sl_2(k)$ and $\rho,\sigma\in\slJ$.

Let $L=M^*\oplus\slJ$ be the Lie superalgebra with product
$$[(m,\rho),(n,\sigma)]=(\omega(\rho,\sigma),[\rho,\sigma]).$$
We have that $L\in\Lst$ since $\omega$ is $sl_2(k)$-invariant. Consider the exact sequence of Lie superalgebras,
$$0\rightarrow M^*\rightarrow L\rightarrow\slJ\rightarrow0.$$
By Theorem \ref {thm:2.3}, the above sequence splits.

Therefore, $H^2(\slJ,M^*)^{sl_2}=0$. $\Box$ 

Here, we summarize the information of $H_r(\slJ)$ for $r=0,1,2$

\begin {pro}\label{pro:3.1}(1) $H_0(\slJ)=k$,

(2) $H_1(\slJ)=k^{D_1|D_2}$, 

(3) $H_2(\slJ)$ is isotypic of type $L(4)$ as an $sl_2(k)$-module.
\end {pro}
{\bf Proof. }(1)It follows from Lemma \ref{lem:3.1}. 

(2)Obviously, $J(D_1|D_2)^2=\oplus_{n\geq 2}J_n(D_1|D_2)$ and $$J(D_1|D_2)/J(D_1|D_2)^2\simeq J_1(D_1|D_2).$$ Note that
$J_1(D_1|D_2)\simeq k^{D_1|D_2}$, where $k^{D_1|D_2}$ is the $k$-superspace with graded dimension $(D_1,D_2)$. By Lemma \ref{lem:3.2},
we obtain (2).

(3)By above Lemma, $$H_2(\slJ)^{sl_2}=H_2(\slJ,L(0))^{sl_2}=0.$$

Let $L(2)$ be the $3$-dimensional simple $sl_2(k)$-module with the trivial action of $\slJ$. Also By above Lemma,
$$H_2(\slJ,L(2))^{sl_2}=0.$$

Since the $\slJ$-actions on $L(2)$ are trivial,  we have that $$H_*(\slJ,L(2))\simeq H_*(\slJ)\otimes L(2)$$ using Chevalley-Eilenberg resolution. 
If $H_2(\slJ)^{\text{ad}}\neq 0$, there exists a submodule $N$ of $H_2(\slJ)^{\text{ad}}$ which is isomorphic to $L(2)$ and
$N\otimes L(2)$ is a submodule of $H_2(\slJ)\otimes L(2)$. However, by Clebsch-Gordan formula,
$$N\otimes L(2)\simeq L(2)\otimes L(2)\simeq L(4)\oplus L(2)\oplus L(0).$$
This means that $H_2(\slJ\otimes L(2)$ has a submodule isomorphic to $L(0)$, which conflicts the fact that $H_2(\slJ,L(2))^{sl_2}=0$. Therefore, we have $H_2(\slJ)^{\text{ad}}=0$.

Also by Chevalley-Eilenberg resolution, the admissible highest weight of $H_2(\slJ)$ as $sl_2(k)$-module are only $0,2,4$. And, we have seen that $0,2$ are impossible. Therefore, $H_2(\slJ)$ is isomorphic to the sum of some $L(4)$. $\Box$  

\

It is a well-known fact that if $\frak m$ is an ordinary free Lie algebra, then $H_r({\frak m})=0$ for $r\geq 2$. Here, $J(D_1|D_2)$ are free objects in the category $\Jsa$ and $\slJ$ are free relative to category $\Lst$ in some sense by Theorem \ref{thm:2.3}. Following \cite{KM}, since only the trivial and adjoint $sl_2$-type occurs in the category $\Lst$, the conjecture seems natural
\begin {conj}\label{con:1}As $sl_2(k)$-modules,
$$H_r(\slJ)^{sl_2}=0 \text{\ \ \  and   }\ \ \ H_r(\slJ)^{\text{ad}}=0$$
for any $r\geq 2$.
\end {conj}

{\bf Remark.} This conjecture is true for $r=2$ by Proposition \ref{pro:3.1}(3).

\section{The graded dimensions of $J_n(D_1|D_2)$}

In this section, we focus on the graded dimension of the homogeneous superspace $J_n(D_1|D_2)$ of $J(D_1|D_2)$. Fixing $D_1,D_2$, for $n\geq 1$, set
$$(\ane,\ano)=(\dim J_n(D_1|D_2)_\bz,\dim J_n(D_1|D_2)_\bo).$$
It is trivial that $(a_1^{\text{even}},a_1^{\text{odd}})=(D_1,D_2)$. What we can say about the vector $(\ane,\ano)$ for $n\geq 2$?

We consider the algebraic group $G=GL(D_1)\times GL(D_2)$, which is the subgroup of the even part of the supergroup $GL(D_1|D_2)$. Then, $J_n(D_1|D_2)$ is a polynomial $G$-module with degree $n$ and $J(D_1|D_2)$ is an analytic $G$-module, described in Section 1 of \cite{KM} or Chapter 1 of \cite{Ma}.  

{\bf Remark.} Here, we only need the even part $G$ of the supergroup $GL(D_1|D_2)$, which is enough for our study on the graded dimension of $J_n(D_1|D_2)$. Moreover, for two $G$-modules $M$ and $N$, there is a natural isomorphism between $G$-modules
$$M\otimes N\rightarrow N\otimes M, m\otimes n\mapsto n\otimes m.$$ 

\

Let $H$ be the subgroup of $G$ constituting by diagonal matrices and $V$ a polynomial $G$-module. For a $(D_1+D_2)$-uple of non-negative integers
$$\text{m}=(m_1,\cdots, m_{D_1},n_1,\cdots,n_{D_2})\in\mZ_{\geq 0}^{D_1+D_2},$$ set
$$V_{\text{m}}=\{v\in V\ |\ h.v=\lambda_1^{m_1}\cdots\lambda_{D_1}^{m_{D_1}}\mu_1^{n_1}\cdots\mu_{D_2}^{n_{D_2}}v,\forall h\in H\},$$
for $h=\text{diag}(\lambda_1,\cdots,\lambda_{D_1},\mu_1,\cdots,\mu_{D_2})\in H$. 

Assume that $V=\oplus_{\text{m}}V_{\text{m}}$ and denote $\deg\text{m}=\sum_{i=1}^{D_1}m_i+\sum_{j=1}^{D_2}n_j$, then
$V=\oplus_{n\geq 0}V_n$ for $V_n=\oplus_{\deg{\text m}=n}V_{\text{m}}$ of degree $n$.

Next, we consider the superspace structure of a $G$-module $V$ of degree $n$. Let $\tilde{Z}$ be the subgroup of $G$ generated by $k^*\mbox{id}$ and $\tau=(\mbox{id}_{D_1},-\mbox{id}_{D_2})$ contained in the center of $G$, which is isomorphic to 
$k^*\times\{\pm 1\}$.

For $\lambda\tau \in\tilde{Z},\lambda\in k^*$, we denote
\begin{align*}
	&V_\bz=\{v\in V\ |\ \lambda \tau.v=\lambda^n v\}=\oplus_{\sum_{j=1}^{D_2}n_j\equiv 0(\mbox{mod }2)}V_{\mbox{m}}\\
	&V_\bo=\{v\in V\ |\ \lambda \tau.v=-\lambda^n v\}=\oplus_{\sum_{j=1}^{D_2}n_j\equiv 1(\mbox{mod }2)}V_{\mbox{m}}.
\end{align*}
We see that $V=V_\bz\oplus V_\bo$.

\subsection{The Grothendieck super-rings ${\mathcal R}^s(K)$ and analytic $G$-modules }

Let $K$ be a reductive algebraic subgroup of $G$, which contains the subgroup $\tilde{Z}$.
For $n\geq 0$, a rational $K$-module on which $\lambda\tau\in\tilde{Z}$ acts by $\pm\lambda^n$ for $\lambda\in k^*$ is called a $K$-module of degree $n$. Denote $\text{Rep}_n(K)$ be the category of $K$-modules of degree $n$.
Let $K_0(\text{Rep}_n(K))$ be the Grothendieck group of the category of homogeneous polynomial $G$-module with degree $n$. Set
\begin{align*}
&{\mathcal R}^s(K)=\sum_{n\geq 0}K_0(\text{Rep}_n(K)),\\
&{\mathcal M}_{>m}^s(K)=\sum_{n>m}K_0(\text{Rep}_n(K)),\\
&{\mathcal M}^s(K)={\mathcal M}_{>0}^s(K).
\end{align*}

Note that there are products $$K_0(\text{Rep}_m(K))\times K_0(\text{Rep}_n(K))\rightarrow K_0(\text{Rep}_{m+n}(K)),$$
induced by the tensor product of the $K$-modules. With the multiplication
$$[(V_\bz,V_\bo)][(W_\bz,W_\bo)]=[(V_\bz\otimes W_\bz\oplus V_\bo\otimes W_\bo,V_\bz\otimes W_\bo\oplus V_\bo\otimes W_\bz)],$$ 
${\mathcal R}^s(K)$ becomes a commutative super-ring($\mZ/2\mZ$-graded ring) and ${\mathcal M}^s(K)$ is a graded ideal.
Moreover, ${\mathcal R}^s(G)$ is complete with respect to the ${\mathcal M}^s(G)$-adic topology, i.e. the topology for which the sequence ${\mathcal M}_{>m}^s(G)$ is a basis of neighborhoods of $0$. Any element $a$ of ${\mathcal R}^s(G)$ can be written as a formal series $a =\sum_{n\geq 0}a_n$ where $a_n\in K_0(\text{Rep}_n(G))$. Note that they are the super versions of ${\mathcal R}(K), {\mathcal M}(K)$ in \cite{KM}.

Furthermore, take $K=\tilde{Z}$. Let $R$ be the ring $\mZ[x]/(x^2-1)$. Identify $(a,b)\in\mZ^2$ with the element $[a+bx]\in R$.
Then, there is an isomorphism of rings $${\mathcal R}^s(\tilde{Z})\rightarrow R[[z]], [V]\mapsto\sum_{n\geq 0}\overline{dim}(V_n)z^n,$$ with ${\mathcal M}^s(\tilde{Z})$ identifying $zR[[z]]$ and ${\mathcal R}^s(\tilde{Z})/{\mathcal M}^s(\tilde{Z})\simeq R$. In the rest of this paper, we always identify the ring ${\mathcal R}^s(\tilde{Z})$ with $R[[z]]$.

A rational $G$-module $V$ is called analytic if $V$ has the decomposition $V=\oplus_{n\geq 0}V_n$ such that $V_n$ is a homogeneous polynomial $G$-module with degree $n$. If $V$ is an analytic $G$-module, then $[V]=\sum_{n\geq 0}[V_n]\in{\mathcal R}^s(G)$. Particular, in the free Jordan superalgebra case,  $J(D_1|D_2), \text{Inn}J(D_1|D_2)$ and ${\mathcal B}^s(J(D_1|D_2))$ are analytic $G$-modules with their natural $\mZ_{>0}$-grading.

\subsection{Effective elements in ${\mathcal R}^s(K)$ and super $\lambda^s$-operation}

The classes $[M]$ of the $K$-modules $M$ are called the effective classes in  ${\mathcal R}^s(K)$. Let ${\mathcal M}^s(K)^+$ be the set of effective classes in  ${\mathcal M}^s(K)$. Then any $a\in {\mathcal R}^s(K)$ can be written as $a_1-a_2$, where $a_1,a_2\in{\mathcal M}^s(K)^+$.

For $a\in{\mathcal M}^s(K)^+$, choose a $K$-module $V=V_\bz\oplus V_\bo$ such that $a=[V]$. Set
$$\lambda^s(a)=\sum_{r\geq 0}(-1)^r\sum_{p+q=r}[\Lambda^pV_\bz\otimes S^qV_\bo]\in1+{\mathcal M}^s(K). $$

It is routine to check that the value $\lambda^s(a)$ is independent of the choice of $V$ and $\lambda^s(a_1+a_2)=\lambda^s(a_1)\lambda^s(a_2)$. Since ${\mathcal R}^s(K)$ is complete, $\lambda^s(a)$ is invertible in  ${\mathcal R}^s(K)$. Then,
\begin{align*}&\lambda([V])=\lambda([V_\bz])\lambda([V_\bo)]\\
=&\left[(\sum_{r\geq 0}(-1)^r\wedge^r V_\bz,0)\right]\left[(\sum_{r\mbox{ even}}S^r V_\bo,-\sum_{r\mbox{ odd}}S^r V_\bo)\right].\end{align*}

Extend $\lambda^s$ on the whole ${\mathcal M}^s(K)$. For $a=a_1-a_2\in{\mathcal M}^s(K)$ with $a_1,a_2\in{\mathcal M}^s(K)^+$, set
$\lambda^s(a)=\lambda^s(a_2)^{-1}\lambda^s(a_1)$. We obtain a group homomorphism $\lambda^s$ from the additive group ${\mathcal M}^s(K)$ to the multiplicative group ${\mathcal R}^s(K)^\times$ of invertible elements in ${\mathcal R}^s(K)$, which is called the super $\lambda^s$-operation on ${\mathcal M}^s(K)$.

\subsection{The main conjecture on structure of $J(D_1|D_2)$}

Let $L(2m)$ be the irreducible $PSL(2)$-module of dimension $2m+1$. Then, $([L(2m)])_{m\geq 0}$ is a $\mZ$-basis of $K_0(PSL(2))$. Since $K_0(K\times PSL(2))=K_0(K)\otimes_\mZ K_0(PSL(2))$, any element $a\in K_0(K\times PSL(2))$ can be written as a sum
$$a=\sum_{m\geq 0}[a:L(2m)][L(2m)],$$
where $[a:L(2m)]\in{\mathcal R}^s(K)$. There are natural generalizations of the concepts of  ${\mathcal R}^s(K\times PSL(2)), {\mathcal M}^s(K\times PSL(2))$ and super $\lambda^s$-operation.

The following Lemma is critical for our purpose.
\begin {lem}\label{lem:4.1} (Main Lemma for super case) Let $K$ be an algebraic subgroup of $G$ which contains $\tilde{Z}$.

(1) There are unique $a^s(K),b^s(K)\in{\mathcal M}^s(K)$ such that the element $$a^s(K)[L(2)]+b^s(K)\in{\mathcal M}^s(K\times PSL(2))$$ satisfies the following two equations :
\begin{align*}
	&[\lambda^s\left(a^s(K)[L(2)]+b^s(K)\right):L(0)]=[k|_K],\\
	&[\lambda^s\left(a^s(K)[L(2)]+b^s(K)\right):L(2)]=-[k^{D_1|D_2}|_K],
\end{align*}
where $k$ is the trivial graded $G$-module, $k^{D_1|D_2}$ is the natural graded $G$-module and $k|_K$ and $k^{D_1|D_2}|_K$ are their restrictions on subgroup $K$.

(2)For $K=G$, set $A(D_1|D_2)=a^s(G)$ and $B(D_1|D_2)=b^s(G)$ in ${\mathcal M}^s(G)$. 

Then, $a^s(K)=A(D_1|D_2)|_K$ and $b^s(G)=B(D_1|D_2)|_K$.
\end {lem}
{\bf Proof. }(1) This follows immediately from the analogous assertion in the Jordan algebra case, which is of a Hensel's Lemma type. See the proof of Lemma 1 in \cite{KM}.

(2) Use the uniqueness of  $a^s(K)$ and $b^s(K)$.   $\Box$

\

Then, we can give a main conjecture for the structure of free Jordan algebra $J(D_1|D_2)$.
\begin {conj}\label{con:2} Let $D_1, D_2$ be non-negative integers with $D_1+D_2\geq 1$. In ${\mathcal R}^s(G)$, we have that
$$[J(D_1|D_2)]=A(D_1|D_2),\ \ [{\mathcal B}^s(J(D_1|D_2))]=B(D_1|D_2),$$
where $A(D_1|D_2)$ and $B(D_1|D_2)$ are defined in the above lemma (2).
\end {conj}

This conjecture is not surprise by the following Theorem.
\begin{thm}\label{thm:4.1} Conjecture 1 implies Conjecture 2.
\end{thm}
{\bf Proof.} Assume that Conjecture 1 holds. 
 
First, note that $\slJ=sl_2(J)\otimes J(D_1|D_2)\oplus{\mathcal B}^s(J(D_1|D_2))$. Hence,
$$[\slJ]=[J(D_1|D_2)][L(2)]+[{\mathcal B}^s(J(D_1|D_2))]\in{\mathcal R}^s(K\times PSL(2)).$$
By the definition of $\lambda^s$, in ${\mathcal R}^s(K)$,
$$\lambda^s([\slJ])=\sum_{r\geq 0}(-1)^r[V_r(\slJ)].$$

Meanwhile, using Euler's characteristic formula for the Grothendieck group of Chevalley-Eilenberg complex, we have
$$\sum_{r\geq 0}(-1)^r[V_r(\slJ)]=\sum_{r\geq 0}(-1)^r[H_r(\slJ)].$$
Hence,
$$\lambda^s([\slJ])=\sum_{r\geq 0}(-1)^r[H_r(\slJ)].$$

In particular, $[\lambda^s([\slJ]):L(0)]=\sum_{r\geq 0}(-1)^r[H_r(\slJ)^{sl_2}]$ and $[\lambda^s([\slJ]):L(2)]=\sum_{r\geq 0}(-1)^r[H_r(\slJ)^{\text{ad}}]$.

By conjecture 1,
\begin{align*}
&[\lambda^s([\slJ]):L(0)]=[H_0(\slJ)^{sl_2}]-[H_0(\slJ)^{sl_2}],\\
&[\lambda^s([\slJ]):L(2)]=[H_0(\slJ)^{\text{ad}}]-[H_1(\slJ)^{\text{ad}}].
\end{align*}

Using Proposition \ref{pro:3.1}(1) and (2), 
\begin{align*}
	&[\lambda^s([\slJ]):L(0)]=[k],\\
	&[\lambda^s([\slJ]):L(2)]=-[k^{D_1|D_2}].
\end{align*}

By the uniqueness of $A(D_1|D_2)$ and $B(D_1|D_2)$ in Lemma \ref{lem:4.1}, we get
$$[J(D_1|D_2)]=A(D_1|D_2),\ \ [{\mathcal B}^s(J(D_1|D_2))]=B(D_1|D_2).\ \ \Box$$

\subsection{The graded dimension of $J_n(D_1|D_2)$}
In this section, we take $K=\tilde{Z}$ and give a more explicit conjecture on the graded dimension of $J_n(D_1|D_2)$.

We have identified the super-ring ${\mathcal R}^s(\tilde{Z})$ with $R[[z]]$. And, $K_0(PSL(2))$ is the subring $\mZ[t+t^{-1}]$ of $\mZ[t,t^{-1}]$ consisting of the symmetric polynomials in $t$ and $t^{-1}$, where $t=e^\alpha$ for the simple root $\alpha$ of the simple algebraic group $PSL(2)$. It follows that
$${\mathcal R}^s(\tilde{Z}\times PLS(2))=R[t+t^{-1}][[z]].$$

For $n>1$, set
$$\overline{\dim}J_n(D_1|D_2)=a^s_n, \ \ \overline{\dim}{\mathcal B}^s_nJ(D_1|D_2)=b^s_n\in R,$$
and the two series $$a^s(z)=\sum_{n\geq 1}a^s_nz^n, \ \ \ b^s(z)=\sum_{n\geq 1}b^s_nz^n\in R[[z]].$$

If Conjecture \ref{con:2} holds, using Lemma \ref{lem:4.1}(2) for $K=\tilde{Z}$, $a^s(z)$ and $b^s(z)$ are uniquely determined by
the two identities in $R[[z]]$:
\begin{align*}
	&[\lambda^s\left(a^s(z)[L(2)]+b^s(z)\right):L(0)]=(1,0),\\
	&[\lambda^s\left(a^s(z)[L(2)]+b^s(z)\right):L(2)]=-(D_1,D_2)z.
\end{align*}

\

The following lemma describes the action of super $\lambda^s$-operation on $R[[z]]$.
\begin {lem}\label{lem:4.2} For $a=(a^{\text{even}},a^{\text{odd}})z^m\in R[[z]]$, we have that
$$\lambda^s(az^m)=((1-z^m)^{a^{\text{even}}},0)\left(\frac{1}{1-z^2},-\frac{z}{1-z^2}\right)^{a^{\text{odd}}}.$$
\end {lem}
{\bf Proof.}  Set $V_\bz(m)$ is the $1$-dimensional $\tilde{Z}$-module such that $\lambda\tau$ acts as $\lambda^m$ and $V_\bo(m)$ is the $1$-dimensional $\tilde{Z}$-module such that $\lambda\tau$ acts as $-\lambda^m$ for $\lambda\in k^*$.

By direct computation,
\begin{align*}
&\lambda^s(az^m)=\lambda^s((z^m,0))^{a^{\text{even}}}\lambda^s((0,z^m))^{a^{\text{odd}}}\\
=&\left([k|_{\tilde{Z}}]-[V_\bz(m)],0\right)^{a^{\text{even}}}\left(\sum_{r\geq0}(-1)^r[S^rV_\bo(m)]\right)^{a^{\text{odd}}}\\
=&\left([k|_{\tilde{Z}}]-[V_\bz(m)],0\right)^{a^{\text{even}}}\left(\sum_{r\text{ even}}[V_\bz(mr)],-\sum_{r\text{ odd}}[V_\bo(mr)]\right)^{a^{\text{odd}}}\\
=&((1-z^m)^{a^{\text{even}}},0)\left(\sum_{r\geq0}z^{2rm},-\sum_{r\geq0}z^{(2r+1)m}\right)^{a^{\text{odd}}}\\
=&((1-z^m)^{a^{\text{even}}},0)\left(\frac{1}{1-z^{2m}},-\frac{z^m}{1-z^{2m}}\right)^{a^{\text{odd}}}. \ \ \ \Box
\end{align*}

We introduce the notation of $t$-integer $[m]_t$ by
$$[m]_t=\frac{t^m-t^{-m}}{t-t^{-1}}=\sum_{i=0}^{m-1}t^{m-1-2i}\in\mZ[t+t^{-1}].$$
Particularly, $[1]_t=1, [2]_t=t^{-1}+t$ and $[3]_t=t^{-2}+1+t^2$. More details of $t$-integer can be found in Chapter 0 of \cite{Ja}.

\begin {lem}\label{lem:4.3}In $R[t+t^{-1}][[z]]$,
\begin{align*}
&\left(\frac{1}{1-t^2z^2},-\frac{tz}{1-t^2z^2}\right)\left(\frac{1}{1-t^{-2}z^2},-\frac{t^{-1}z}{1-t^{-2}z^2}\right)\\
=&\left(\sum_{i=0}^\infty[2i+1]_tz^{2i},-\sum_{i=0}^\infty[2i+2]_tz^{2i+1}\right).
\end{align*}
\end {lem}
{\bf Proof. } By the identity $$\frac{1}{(1-t^2z^2)(1-t^{-2}z^2)}=\sum_{i=0}^\infty\left(\sum_{j=0}^it^{2(i-2j)}\right)z^{2i},$$
we have
\begin{align*}
&\frac{1+z^2}{(1-t^2z^2)(1-t^{-2}z^2)}=\sum_{i=0}^\infty\left(\sum_{j=0}^it^{2(i-j)}\right)z^{2i}=\sum_{i=0}^\infty[2i+1]_tz^{2i},\\
&\frac{z(t+t^{-1})}{(1-t^2z^2)(1-t^{-2}z^2)}=\sum_{i=0}^\infty\left(\sum_{j=0}^it^{2(i-j)+1}\right)z^{2i}=\sum_{i=0}^\infty[2i+2]_tz^{2i+1}.
\end{align*}
The result is direct. $\Box$

For two series $a^s(z)=\sum_{n\geq 1}(a_n^{even},a_n^{odd})z^n$, $b^s(z)=\sum_{n\geq 1}(b_n^{even},b_n^{odd})z^n$ in $zR[[z]]$, define the series
\begin{align*}
	&\Phi^s=\prod_{n\geq 1}[\left((1-[2]_tz^n+z^{2n})^{a_n^{even}}(1-z^n)^{a_n^{even}+b_n^{even}},0\right)\cdot\\
	&\left(\sum_{i=0}^\infty[2i+1]_tz^{2in},-\sum_{i=0}^\infty[2i+2]_tz^{(2i+1)n}\right)^{a_n^{odd}}\left(\frac{1}{1-z^{2n}},-\frac{z^n}{1-z^{2n}}\right)^{a_n^{odd}+b_n^{odd}}]\\
	&\in R[t+t^{-1}][[z]].
\end{align*}

\begin{lem}\label{lem:4.4}There are uniquely two series $a^s(z)$ and $b^s(z)\in zR[[z]]$  such that the associated $\Phi$ satisfies the following two equations:
	\begin{align*}
		&\text{Res}_{t=0}(t^{-1}-1)\Phi^s dt=(1,0),\\
		&\text{Res}_{t=0}(1-t)\Phi^s dt=-(D_1,D_2)z.
	\end{align*}
\end {lem}

{\bf Proof.} First, we see that for $a=\sum_i a_it^i\in\mZ[t+t^{-1}]=K_0(PSL(2))$, $[a:L(0)]=a_0-a_{-1}$ and $[a:L(2)]=a_{-1}-a_{-2}$. Thus, $[a:L(0)]=\text{Res}_{t=0}(t^{-1}-1)a$ and $[a:L(2)]=\text{Res}_{t=0}(1-t)a$.

In Lemma \ref{lem:4.1}(2), taking $K=\tilde{Z}$, we only need to show that 
$\Phi^s=\lambda^s\left(a^s(z)[L(2)]+b^s(z)\right)$. See that $[L(0)]=1$ and $[L(2)]=t^{-1}+1+t$. By Lemma \ref{lem:4.2} and \ref{lem:4.3},
\begin{align*}
	&\lambda^s\left(a^s(z)[L(2)]+b^s(z)\right)\\
	=&\lambda^s(\sum_{n\geq 1}(a_n^{even},a_n^{odd})(t+1+t^{-1})z^n+\sum_{n\geq 1}(b_n^{even},b_n^{odd})z^n)\\
	=&\prod_{n\geq 1}\lambda^s((a_n^{even},a_n^{odd})(t+1+t^{-1})z^n)\cdot\lambda^s((b_n^{even},b_n^{odd})z^n)\\
	=&\Phi^s. \ \ \ \ \Box
\end{align*}

We give
\begin {conj}\label{con:3} Let $a^s(z)=\sum_{n\geq 1}a^s_nz^n$, $b^s(z)=\sum_{n\geq 1}b^s_nz^n\in R[[z]]$ be two series defined in Lemma \ref{lem:4.4}.
Then, $a^s(z)$ and $b^s(z)$ are given by the graded dimensions of $J_n(D_1|D_2)$ and ${\mathcal B}_n^s(J(D_1|D_2))$.
\end {conj}

\begin{thm}\label{thm:4.2} Conjecture 2 implies Conjecture 3.
\end{thm}
{\bf Proof.} Assuming the Conjecture 2 holds, we have that
$$[J(D_1|D_2)]=A(D_1|D_2),\ \ [{\mathcal B}^s(J(D_1|D_2))]=B(D_1|D_2).$$
Restricting on $\tilde{Z}$, $a^s(z)=A(D_1|D_2)|_{\tilde{Z}},  b^s(z)=B(D_1|D_2)|_{\tilde{Z}}$
are satisfied the two equations in Lemma \ref{lem:4.3}. $\Box$

\

Furthermore, we can take a linear combination of the two equations and cancel the $b^s(z)$. Set the series $\Psi^s$ in $R[t+t^{-1}][[z]]$
\begin{align*}
	&\Psi^s=\prod_{n\geq 1}\left(1-[2]_tz^n+z^{2n},0\right)^{a_n^{even}}\\
	&\left(\sum_{i=0}^\infty[2i+1]_tz^{2in},-\sum_{i=0}^\infty[2i+2]_tz^{(2i+1)n}\right)^{a_n^{odd}},
\end{align*}
and the series $\Theta^s$ in $R[[z]]$
$$\Theta^s=\prod_{n\geq 1}\left(1-z^n,0\right)^{a_n^{even}+b_n^{even}}\left(\frac{1}{1-z^{2n}},-\frac{z^n}{1-z^{2n}}\right)^{a_n^{odd}+b_n^{odd}}
.$$
Then, $\Phi^s=\Theta^s\Psi^s$. Multiplying the first equation by $(D_1,D_2)z$ and adding it on the second one, the series
$a^s(z)$ satisfies the equation in $R[[z]]$
$$\Theta^s(\text{Res}_{t=0}((D_1z,D_2z)(t^{-1}-1)+(1-t,0))\Psi^s dt)=0.$$
Since $R[[z]]$ is an integral domain and $\Theta^s\equiv 1(\mod z)$ is nonzero, we can cancel it and obtain
$$\text{Res}_{t=0}\psi^s\Psi^s dt=0,$$
where $\psi^s=(D_1z,D_2z)t^{-1}+(1-D_1z,-D_2z)+(-1,0)t\in R[t^{\pm1}][z]$.

\begin{lem}\label{lem:4.5}There exists a series $a^s(z)\in zR[[z]]$ determined uniquely by the following equations:
$$\text{Res}_{t=0}\psi^s\Psi^s dt=0$$
for above $\Psi^s$ and $\psi^s$.\ \ \ \ $\Box$
\end {lem}

Naturally,
\begin {conj}\label{con:4} Set $a_n^s=\overline{\dim}J_n(D_1|D_2)$. The sequence $a^s_n=(\ane,\ano)$ is the unique solution of the following equation:
\begin{equation}\text{Res}_{t=0}\psi^s\Psi^s dt=0.\tag{{$\mathcal E$}}\end{equation}
\end {conj}

\section{Some evidences for these Conjectures}

In this section, we consider the free Jordan superalgebra $J(D_1|D_2)$ for some specialized $(D_1,D_2)$. 

First, for $D_2=0$, $J(D_1|0)$ becomes the Jordan algebra $J(D_1)$ and $sl_2(J(D_1|0))$ is the $\TAG$ Lie algebra $sl_2J(D_1)$. Hence, $\ano=b_n^{\text{odd}}=0$ for all $n\geq 1$. In these cases, our Conjecture \ref{con:1}-\ref{con:4} are nothing but the Conjecture 2 and 1(including weak and weakest version) in \cite{KM}.

\subsection{The cases $D_1+D_2=1$}

For the case $(D_1,D_2)=(0,1)$, $J(0|1)=k.y_1$ is the $1$-dimensional zero algebra over $k$, namely $J(0|1)^2=0$. Moreover, ${\mathcal B}^s(J(0|1))=k.\{y_1\otimes y_1\}$.
Then, the only nonzero coefficients of $a^s(z)$ and $b^s(z)$ are $a^s_1=(0,1),b^s_2=(1,0)$. And, $sl_2(J(0|1))$ is a $(1,3)$-dimensional Lie superalgebra. Now,
\begin{align*}
\Phi^s=&(1-z^2,0)\left(\sum_{i=0}^\infty[2i+1]_tz^{2i},-\sum_{i=0}^\infty[2i+2]_tz^{2i+1}\right)\left(\frac{1}{1-z^2},-\frac{z}{1-z^2}\right)\\
=&(1,-z)\left(\sum_{i=0}^\infty[2i+1]_tz^{2i},-\sum_{i=0}^\infty[2i+2]_tz^{2i+1}\right).
\end{align*}

Since $\text{Res}_{t=0}t^{-1}[m]_t=\begin{cases}&0, m\equiv0(\mod 2)\\&1, m\equiv1(\mod 2)\end{cases}$, and $ \text{Res}_{t=0}[m]_t=\begin{cases}&1, m\equiv0(\mod 2)\\&0, m\equiv1(\mod 2)\end{cases}$,
$$\text{Res}_{t=0}(t^{-1}-1)\Phi^s dt=(1,-z)\left(\frac{1}{1-z^2},\frac{z}{1-z^2}\right)=(1,0)$$
and similarly,
$$\text{Res}_{t=0}(1-t)\Phi^s dt=(1,-z)\left(\frac{-z^2}{1-z^2},\frac{-z}{1-z^2}\right)=-(0,1)z.$$
Hence, Conjecture 3 holds for $J(0|1)$.

Next, we check Conjecture 4 directly as follows. Now,
$$\psi^s=(0,z)t^{-1}+(1,-z)+(-1,0)t$$ and $$\Psi^s=\left(\sum_{i=0}^\infty[2i+1]_tz^{2i},-\sum_{i=0}^\infty[2i+2]_tz^{2i+1}\right)$$  Set 
$\Psi^s=\sum_i\Psi^s_it^i$. We have that 
$\Psi^s_0=(\frac{1}{1-z^2},0), \Psi^s_{-1}=(0,-\frac{z}{1-z^2})$ and 
$\Psi^s_{-2}=(\frac{z^2}{1-z^2},0)$.
Then,
$$\text{Res}_{t=0}\psi^s\Psi^s dt =(0,z)\Psi^s_0+(1,-z)\Psi^s_{-1}+(-1,0)\Psi^s_{-2}=0.$$

\begin{pro}\label{prop:5.1} Conjecture 3 and 4 hold for the free Jordan superalgebra with only one generator.
\end{pro}
{\bf Proof. }If the generator is even, the Jordan algebra case can be found in Section 2.6 of \cite{KM}. 
And, if the generator is odd, the case has been discussed as above. $\Box$

\begin{rem}\label{rem:4.1} Here, we point out a different phenomenon of free Jordan superalgebras from the free Jordan algebras.
For a free Jordan algebra $J(D)$, the Conjecture 2 on the homology of $sl_2J(D)$ implies that
$\text{Inn}(J(D))={\mathcal B}(J(D))$(See Lemma 12 in \cite{KM}). But for Jordan superalgebra $J(0|1)$, we have seen that
${\mathcal B}^s(J(0|1))=k.\{y_1\otimes y_1\}\neq 0$, but $\text{Inn}(J(0|1))=0$ since $J(0|1)$ is just a zero algebra.
\end{rem}

\subsection{The cases $D_1+D_2=2$}

The structures of free Jordan superalgebras with two generators are more complex. We have only incomplete results.

For the Jordan algebra case $(D_1,D_2)=(2,0)$, it is shown that $\text{Res}_{t=0}\psi\Psi dt\equiv0(\mod z^{16})$ to support Conjecture 4 in Section 2.6 of \cite{KM}.

For the case $(D_1,D_2)=(0,2)$, we have
 \begin{align*}J_1(0|2)=&k.\{y_1,y_2\},\\
 	J_2(0|2)=&k.\{y_1\cdot y_2\},\\
 	J_3(0|2)=&k.\{(y_1\cdot y_2)\cdot y_1,(y_1\cdot y_2)\cdot y_2\},\\
 	J_4(0|2)=&k.\{(y_1\cdot y_2)\cdot(y_1\cdot y_2),((y_1\cdot y_2)\cdot y_1)\cdot y_1,((y_1\cdot y_2)\cdot y_1)\cdot y_2,\\
 	&((y_1\cdot y_2)\cdot y_2)\cdot y_1,((y_1\cdot y_2)\cdot y_2)\cdot y_2\}.
 \end{align*}
Since the super Jordan identities are all trivial up to $n\leq 4$ in this case, we get
$$a^s_1=(0,2),\ a^s_2=(1,0),\ a^s_3=(0,2),\ a^s_4=(5,0).$$

We can compute $\Psi^s(\mod z^5)$ by hand, using
$[2]_t^2=[3]_t+1,[2][3]=[4]_t+[2]_t,[3]^2_t=[5]_t+[3]_t+1,[2]_t[4]_t=[5]_t+[3]_t$
and $[2]_t^3=[4]_t+2[2]_t$.
\begin{align*}\Psi^s\equiv&((1-[2]_tz^2+z^4)(1-5[2]_tz^4),0)\cdot\\
	&(1+[3]_tz^2+[5]_tz^4,-[2]_tz-[4]_tz^3)^2(1,-[2]_tz^3)^2\\
\equiv&(1-[2]_tz^2+(1-5[2]_t)z^4,-2[2]_tz^3)\cdot\\
&(1+([2]_t^2+2[3]_t)z^2+([3]^2_t+2[5]_t+2[2]_t[4]_t)z^4,\\
&-2([2]_tz+([2]_t[3]_t+[4]_t)z^3))\\
\equiv&(1+([2]_t^2+2[3]_t-[2]_t)z^2+\\
&(1-5[2]_t+[3]^2_t+2[5]_t+2[2]_t[4]_t-[2]_t^3-2[2]_t[3]_t+4[2]^2_t)z^4,\\
&-2([2]_tz+([2]_t[3]_t+[4]_t-[2]_t^2+[2]_t)z^3)\\
\equiv&((1+(1-[2]_t+3[3]_t)z^2+(6-9[2]_t+7[3]_t-3[4]_t+5[5]_t)z^4,\\
&-2([2]_tz+(-1+2[2]_t-[3]_t+2[4]_t))z^3)(\mod z^5).
\end{align*}
Then,
 \begin{align*}&\Psi^s_0\equiv(1+4z^2+18z^4,4z^3)(\mod z^5),\\
 	&\Psi^s_{-1}\equiv(-z^2-12z^4,-2z-8z^3)(\mod z^5),\\
 	&\Psi^s_{-2}\equiv(3z^2+12z^4,2z^3)(\mod z^5).
\end{align*}
Furthermore, using
$$\psi^s=(0,2z)t^{-1}+(1,-2z)+(-1,0)t,$$
we have $$\text{Res}_{t=0}(\psi^s\Psi^s) dt=(0,2z)\Psi^s_0+(1,-2z)\Psi^s_{-1}+(-1,0)\Psi^s_{-2}\equiv0(\mod z^5).$$

\

For $(D_1,D_2)=(1,1)$, it is obvious that
\begin{align*}J_1(1|1)=&k.\{x\}\oplus k.\{y\},\\
	J_2(1|1)=&k.\{x^2\}\oplus k.\{x\cdot y\},\\
	J_3(1|1)=&k.\{x^3,(x\cdot y)\cdot y\}\oplus k.\{x^2\cdot y,(x\cdot y)\cdot x\}.
\end{align*}
And, $a^s_1=(1,1),\ a^s_2=(1,1),\ a^s_3=(2,2)$.
The super Jordan identities are effective for degree $4$, which deduce two equations,
\begin{align*}&[L_{x^2},L_x]=0,\\
	&[L_{x^2},L_y]+2[L_{x\cdot y},L_x]=0.\\
\end{align*}
Acting on $x$ and $y$, $(x\cdot y)\cdot x^2=(x^2\cdot y)\cdot x$, $3(x\cdot y)\cdot x^3-x^3\cdot y-2((x\cdot y)\cdot x)\cdot x=0$ and $(x^2\cdot y)\cdot y+2((x\cdot y)\cdot y)\cdot x=0$.
Hence, $a^s_4=(3,3)$.

To compute $\Psi^s(\mod z^5)$, one has
\begin{align*}\Psi^s\equiv&((1-[2]_tz+z^2)(1-[2]_tz^2+z^4)(1-2[2]_tz^3)(1-3[2]_tz^4),0)\cdot\\
	&(1+[3]_tz^2+[5]_tz^4,-[2]_tz-[4]_tz^3)(1+[3]_tz^4,-[2]_tz^2)(1,-2[2]_tz^3)(1,-3[2]_tz^4)\\
\equiv&(1-[2]_tz+(1-[2]_t)z^2+([2]_t^2-2[2]_t)z^3+(1-4[2]_t+2[2]_t^2)z^4,0)\\
&(1+[3]_tz^2+[2]_t^2z^3+(2[2]_t^2+[3]_t+[5]_t)z^4,\\
&-([2]_tz+[2]_tz^2+(2[2]_t+[4]_t)z^3+([2]_t[3]_t+3[2]_t)z^4))\\
\equiv&(1-[2]_tz+(1-[2]_t+[3]_t)z^2+(2-3[2]_t+2[3]_t-[4]_t)z^3\\
&+(5-7[2]_t+6[3]_t-2[4]_t+[5]_t)z^4,-[2]_tz+(1-[2]_t+[3]_t)z^2+\\
&(2-3[2]_t+2[3]_t-[4]_t)z^3+(5-7[2]_t+6[3]_t-2[4]_t+[5]_t)z^4)(\mod z^5).
\end{align*}
and
 \begin{align*}&\Psi^s_0\equiv(1+2z^2+4z^3+12z^4,2z^2+4z^3+12z^4)(\mod z^5),\\
	&\Psi^s_{-1}\equiv-(z+z^2+4z^3+9z^4,z+z^2+4z^3+9z^4)(\mod z^5),\\
	&\Psi^s_{-2}\equiv(z^2+2z^3+7z^4,z^2+2z^3+7z^4)(\mod z^5).
\end{align*}
Note that 
$\psi^s=(z,z)t^{-1}+(1-z,-z)+(-1,0)t$. Finally, we have
$$\text{Res}_{t=0}(\psi^s\Psi^s) dt=(z,z)\Psi^s_0+(1-z,-z)\Psi^s_{-1}+(-1,0)\Psi^s_{-2}\equiv0(\mod z^5).$$
We obtain the following proposition, which gives some evidences of conjectures \ref{con:4} when $D_1+D_2=2$.
\begin{pro}\label{prop:5.2} For $\Psi^s,\psi^s$ given above, $D_1+D_2=2$,
$$\text{Res}_{t=0}\psi^s\Psi^s dt\equiv0\begin{cases}(\mod z^{16}), &(D_1,D_2)=(2,0)\\(\mod z^5), &(D_1,D_2)=(1,1)\\(\mod z^5), &(D_1,D_2)=(0,2)\end{cases}.\Box$$
\end{pro}

\section{Concluding Remarks}
In the Section 2, we have proved a main result on the adjoint pair of functor between category $\Jsa$ and $\Lst$.
\begin{thm}\label{thm:6.1} There are two functors $\TAG:\Jsa\rightarrow\Lst$ and Tits functor $T:\Lst\rightarrow\Jsa$ which form a pair of adjoint functors,  namely,
	$$\Hom_\Lst(\TAG(J),{\mathfrak g})=\Hom_\Jsa(J,T({\mathfrak g}))$$
	for any $J\in\Jsa$ and ${\mathfrak g}\in\Lst$.
\end{thm}

For the free Jordan superalgebra $J(D_1|D_2)$ and its $\TAG$ Lie superalgerba $\slJ$, there are $4$ conjectures listed as follows,
\begin {conje}(The conjecture on the homology of $\slJ$) 
As $sl_2(k)$-modules, the homology superspaces of $\slJ$ satisfy
$$H_r(\slJ)^{sl_2}=0 \text{\ \ \  and   }\ \ \ H_r(\slJ)^{\text{ad}}=0$$
for any $r\geq 2$.
\end {conje}

\begin {conje}(The conjecture on the structure of $J(D_1|D_2)$) 

Let $J(D_1|D_2)$ be free Jordan superalgebra with $D_1+D_2\geq 1$. In ${\mathcal R}^s(G)$, set
$$A(D_1|D_2)=[J(D_1|D_2)],\ \ B(D_1|D_2)=[{\mathcal B}^s(J(D_1|D_2))].$$
Then, $A(D_1|D_2)$ and $B(D_1|D_2)$ are uniquely determined by two equations in ${\mathcal R}^s(G)$:
\begin{align*}
	&[\lambda^s\left(A(D_1|D_2)L(2)+B(D_1|D_2)\right):L(0)]=[k],\\
	&[\lambda^s\left(A(D_1|D_2)L(2)+B(D_1|D_2)\right):L(2)]=-[k^{D_1|D_2}],
\end{align*}
where $k$ is seen as the trivial graded $G$-module, $k^{D_1|D_2}$ is the natural graded $G$-module.
\end {conje}

\begin {conje}(The conjecture on the graded dimensions of both $J_n(D_1|D_2)$ and ${\mathcal B}^s_n(J_n(D_1|D_2))$) 

Set $R=\mZ[x]/(x^2-1)$. Let $a_n^s=\overline{\dim}J_n(D_1|D_2)$, $b^s_n=\overline{\dim}{\mathcal B}^s_n(J(D_1|D_2))$ be two sequences in $R$.
Then, the series $a^s(z)=\sum_{n\geq 1}a_n^s z^n$ and $b^s(z)=\sum_{n\geq 1}b_n^s z^n$ are uniquely determined by two identities in $R[[z]]$:
\begin{align*}
	&\text{Res}_{t=0}(t^{-1}-1)\Phi^s dt=(1,0),\\
	&\text{Res}_{t=0}(1-t)\Phi^s dt=-(D_1,D_2)z,
\end{align*}
where $\Phi^s\in R[t+t^{-1}][[z]]$ is given by
\begin{align*}
	&\Phi^s=\prod_{n\geq 1}[\left((1-[2]_tz^n+z^{2n})^{a_n^{even}}(1-z^n)^{a_n^{even}+b_n^{even}},0\right)\cdot\\
	&\left(\sum_{i=0}^\infty[2i+1]_tz^{2in},-\sum_{i=0}^\infty[2i+2]_tz^{(2i+1)n}\right)^{a_n^{odd}}\left(\frac{1}{1-z^{2n}},-\frac{z^n}{1-z^{2n}}\right)^{a_n^{odd}+b_n^{odd}}]
\end{align*}
\end {conje}

\begin {conje}(The conjecture on the graded dimension of $J_n(D_1|D_2)$) 

Set $a_n^s=\overline{\dim}J_n(D_1|D_2)$. The sequence $a^s_n=(\ane,\ano)$ is the unique solution of the following equation $(\mathcal E)$:
\begin{equation}\text{Res}_{t=0}\psi^s\Psi^s dt=0,\tag{{$\mathcal E$}}\end{equation}
where
\begin{align*}
	&\Psi^s=\prod_{n\geq 1}\left(1-[2]_tz^n+z^{2n},0\right)^{a_n^{even}}\\
	&\left(\sum_{i=0}^\infty[2i+1]_tz^{2in},-\sum_{i=0}^\infty[2i+2]_tz^{(2i+1)n}\right)^{a_n^{odd}}
\end{align*}
and $\psi^s=(D_1z,D_2z)t^{-1}+(1-D_1z,-D_2z)+(-1,0)t\in R[t^{\pm1}][z]$.
\end {conje}

In this paper, we have proven that Conjecture 1 $\Longrightarrow$ 2 $\Longrightarrow$ 3 $\Longrightarrow$ 4.
Conjecture 3 and 4 are true when $D_1+D_2=1$ and some evidences are given when $D_1+D_2=2$.

\

\end{document}